\documentclass{article}
\usepackage{amsmath}
\usepackage{amsmath,amssymb,latexsym,color}
\usepackage[mathscr]{eucal}
\newtheorem{thm}{Theorem}[section]
\newtheorem{prop}[thm]{Proposition}
\newtheorem{lemma}[thm]{Lemma}
\newtheorem{cor}[thm]{Corollary}
\newtheorem{defin}[thm]{Definition}

\newtheorem{example}{Example}[section]

\newcommand{\proof}{{\it Proof.\quad}}
\newcommand{\qed}{\hfill\Box\medskip}

\usepackage{CJK}
\begin{document}

\renewcommand{\baselinestretch}{1.2}

\title{\bf Vector spaces over finite commutative rings}

\author{
Jun Guo\thanks{Corresponding author. guojun$_-$lf@163.com}\quad Junli Liu\thanks{Junli810@163.com}\quad Qiuli Xu\thanks{qiulixu1008@126.com}\\
{\footnotesize Department of Mathematics, Langfang Normal University, Langfang  065000,  China} }
 \date{ }
 \maketitle

\begin{abstract}
 Vector spaces over finite fields and Anzahl formulas of subspaces were studied by Wan (Geometry of Classical Groups over Finite Fields, Science Press, 2002). As a generalization, we study  vector spaces and singular linear spaces
 over commutative rings, and obtain some  Anzahl formulas and dimensional formula for subspaces. Moreover, we discuss
 arcs and caps by using these subspaces.



\end{abstract}

\section{Introduction}
Throughout the paper, our rings always are finite and contain the identity $1\not=0$.

Let $R$ be a commutative ring and $R^\ast$ denote its unit group.
For a subset $S$ of $R$, let $S^{m\times n}$ be the set of all $m\times n$ matrices over $S$, and $S^{n}=S^{1\times n}$.
A matrix in $S^n$ is also called an $n$-dimensional row vector  over $S$.
 Let $I_r$ ($I$ for short) be the $r\times r$ identity matrix, $0_{m,n}$ ($0$ for short) the $m\times n$ zero matrix and $0_n=0_{n,n}$,
  and $\det(X)$ the determinant of a square matrix $X$ over $R$. The set of $n\times n$ invertible matrices forms a group under matrix multiplication, called the {\it general linear group} of degree $n$ over $R$ and denoted by $GL_n(R)$.

  Let $A\in R^{m\times n}$. Denote by $I_k(A)$ the ideal in $R$
generated by all $k\times k$ minors of $A$.
Let $\hbox{Ann}_{R}(I_k(A)) =\{x\in R:xa=0\;\hbox{for all}\;a\in I_k(A)\}$ denote the annihilator of $I_k(A)$. The {\it McCoy rank} of $A$, denoted by $\hbox{rk}(A)$, is the following integer:
$$\hbox{rk}(A) = \max\{k : \hbox{Ann}_{R}(I_k(A)) =\{0\}\}.$$
For $m\leq n$, let $R^{m\times n}(m)$ denote the set of all $m\times n$ matrices of McCoy rank $m$ over $R$.
Then $GL_{n}(R)=R^{n\times n}(n)$. For $0\leq m_1\leq m\leq n$ and a fixed  $A_1\in R^{m_1\times n}(m_1)$, we define
 \begin{equation}\label{equa-zj}
 R^{m\times n}(m,A_1)=\left\{\left(\begin{array}{c}A_1\\ A_2\end{array}\right)\in R^{m\times n}(m)
 : A_2\in R^{(m-m_1)\times n}(m-m_1)\right\}.
 \end{equation}

\begin{prop}\label{prop1.1}{\rm(See \cite{Brown}).}
Let $A\in R^{m\times n}$. Then
\begin{itemize}
\item[\rm(i)]
$0\leq {\rm rk}(A) \leq \min\{m, n\}$.

\item[\rm(ii)]
${\rm rk}(A) = {\rm rk}(PAQ)$ for all $P\in GL_m(R)$ and $Q\in GL_n(R)$.

\item[\rm(iii)]
If $m\leq n$, then ${\rm rk}(A)=m$ if and only if the set of all the rows of $A$ is linearly independent.
\end{itemize}
\end{prop}

  For $A\in R^{m\times n}$ and $B\in R^{n\times m}$, if $AB=I_m$, we say that $A$ has
a {\it right inverse} and $B$ is a right inverse of $A$.
The free module $R^{n}$ of rank $n$ is called the $n$-dimensional {\it row vector space} over $R$. Let $\alpha_i\in R^{n}$ for each $i\in[m]:=\{1,\ldots,m\}$.
The vector subset $\{\alpha_1,\alpha_2,\ldots,\alpha_m\}$ is called {\it unimodular} if the matrix
$(\alpha_1^{\top}, \alpha_2^{\top}, \ldots, \alpha_m^{\top})^{\top}$
has a right inverse, where $X^{\top}$ is the transpose of $X$.
Note that $\{\alpha\}$ is unimodular, where $\alpha=(r_1,\ldots,r_n)\in R^n$, if and only if
the ideal generated by $r_1,\ldots,r_n$ is equal to $R$.
Let $V\subseteq R^{n}$ be a {\it linear subset} (i.e. $R$-module). A {\it largest unimodular vector subset}
of $V$ is a unimodular vector subset of $V$ which has maximum number of vectors. The
{\it dimension} of $V$, denoted by $\dim(V)$, is the number of vectors in a largest unimodular
vector subset of $V$. Clearly, $V$ contains a subspace of dimension
$\dim(V)$ and $\dim(V)=0$ if and only if $V$ does not contain a unimodular
vector. If a linear subset $X$ of $R^{n}$  has a unimodular basis with $m$ vectors, then $X$ is
called an $m$-{\it dimensional vector subspace} ($m$-{\it subspace} for short) of $R^{n}$. We define the $0$-subspace to be $\{0_{1,n}\}$.

Let $\mathbb{F}_q$ be a finite field with $q$ elements. Let $\mathbb{F}_q^{n}$ be one of the $n$-dimensional classical
spaces  over $\mathbb{F}_q$, and let $G_{n}$ be the corresponding classical group of degree $n$, see Wan's monograph \cite{Wan2}.
Wan  determined the orbits of subspaces of $\mathbb{F}_q^{n}$ under the action of $G_n$,  discussed the following questions
and found their the answers:
\begin{itemize}
\item[\rm(i)]
How should the orbits be described?

\item[\rm(ii)]
What are the lengths of the orbits?

\item[\rm(iii)]
What is the number of subspaces in any orbit contained in a given subspace?

\item[\rm(iv)]
What is the number of subspaces in any orbit containing a given subspace?
\end{itemize}
As generalizations, Wang, the first present author and Li \cite{Wang,Wang2}
studied problems (i)-(iv) and found the answers in the singular linear space over  $\mathbb{F}_{q}$;
Huang, Lv and Wang \cite{Huang3} studied problems (ii)-(iv) and found the answers in the vector space over the residue
class ring $\mathbb{Z}_{p^s}$, where $p$ is a prime number;
The first present author \cite{Guo3}  studied problems (ii)-(iv) and found the  answers in the vector space
over the residue class ring $\mathbb{Z}_{m}$;
Sirisuk and Meemark \cite{Sirisuk} studied problem (ii) by using  a lifting technique,
and found the answer in the vector space over the commutative ring $R$.

In this paper, we study problems (i)-(iv) in the vector space and the singular linear space
 over the commutative ring $R$
and find all the answers by using a matrix method. The structure is as follows.
 In Sections~2 and~3, we obtain some Anzahl formulas of matrices over $R$ for later references.
 In Section~4, we find the answers of problems (i)-(iv) in the vector space, and obtain dimensional formula for subspaces.
 In Section~5, we find the answers of problems (i)-(iv) in the singular linear space.
 As an application, we discuss arcs and caps  in Section~6.

\section{Matrices over local rings}
A {\it local ring} is a commutative ring which has a unique maximal ideal.  For a local ring $R$,
its unique maximal ideal is $M = R\setminus R^\ast$,
and the residue field of $R$ is the quotient ring $R/M$. In this section, we always assume that
$R$ is a local ring with the unique maximal ideal $M$
and the residue field $R/M$.

Let  $\pi$ be the natural surjective homomorphism from
$(R,+,\cdot)$  to $(R/M,+,\cdot)$, i.e.
$$\begin{array}{rcl}
\pi:R &\longrightarrow& R/M,\\
 r &\longmapsto& r+M.
\end{array}$$
For $A=(a_{uv})\in R^{m\times n}$, let
$\pi(A)=(\pi(a_{uv}))$.
For matrices $A\in R^{m\times n}$ and $B\in R^{n\times k}$, it is easy to see that
 \begin{equation}\label{shi1}
 \pi(AB)=\pi(A)\pi(B).
 \end{equation}

 \begin{lemma}\label{lemma2.1}{\rm(See \cite[Lemma~2]{Brawley}).}
 Let $A\in R^{m\times n}$. Then ${\rm rk}(A)={\rm rk}(\pi(A))$.
 \end{lemma}

 \begin{cor}\label{lemma2.1-2}
 Let $s\leq n$, and $x_1,x_2,\ldots,x_s\in R^n$. Then $\{x_1,x_2,\ldots,x_s\}$ is linearly independent
over $R$ if and only if $\{\pi(x_1),\pi(x_2),\ldots,\pi(x_s)\}$ is linearly independent over $R/M$.
 \end{cor}
 \proof Write $A=(x_1^\top,x_2^\top,\ldots,x_s^\top)^\top$. By Proposition~\ref{prop1.1} and Lemma~\ref{lemma2.1}, $\{x_1,x_2,\ldots,x_s\}$ is linearly independent over $R$ if and only if ${\rm rk}(A)={\rm rk}(\pi(A))=s$, if and only if $\{\pi(x_1),\pi(x_2),$ $\ldots,\pi(x_s)\}$ is linearly independent over $R/M$. $\qed$

Recall that $\{\alpha\}$ is unimodular, where $\alpha=(r_1,\ldots,r_n)\in R^n$, if and only if
the ideal generated by $r_1,\ldots,r_n$ is equal to $R$. So, we obtain the following result.

\begin{lemma}\label{lemma2.2}{\rm (See \cite[Proposition~1.1]{Sirisuk}).}
Let $x=(r_1,\ldots,r_n)\in R^n$. Then the following statements are equivalent.
\begin{itemize}
\item[\rm(i)]
$\{x\}$ is unimodular.

\item[\rm(ii)]
$r_i$ is a unit for some $i\in \{1,\ldots, n\}$.

\item[\rm(iii)]
$\{x\}$ is  linearly independent.
\end{itemize}
\end{lemma}

\begin{lemma}\label{lemma2.3}
Let $m\leq n$ and $A\in R^{m\times n}$. Then ${\rm rk}(A)=m$ if and only if there exists an $S\in GL_n(R)$ such that
$AS=(I_m,0_{m,n-m})$.
\end{lemma}
\proof If  there exists an $S\in GL_n(R)$ such that
$AS=(I_m,0_{m,n-m})$, by Proposition~\ref{prop1.1} (ii), ${\rm rk}(A)={\rm rk}(I_m,0_{m,n-m})$.
From definition of McCoy rank of matrices, we deduce that ${\rm rk}(I_m,0_{m,n-m})\geq m$.
By Proposition~\ref{prop1.1} (i), ${\rm rk}(I_m,0_{m,n-m})\leq m$. It follows that ${\rm rk}(A)={\rm rk}(I_m,0_{m,n-m})=m$.

Suppose ${\rm rk}(A)=m$. Let $A=(a_{ij})=(\alpha_1^{\top},\alpha_2^{\top},\ldots,\alpha_m^{\top})^{\top}$. By Proposition~\ref{prop1.1} (iii), $\{\alpha_1,\alpha_2,\ldots,\alpha_m\}$
  is linearly independent, which implies that $\{\alpha_1=(a_{11},a_{12},\ldots,a_{1n})\}$ is linearly independent. Without loss of generality, by Lemma~\ref{lemma2.2},
we may assume that $a_{11}$ is a unit. Then
$$
S_1=\left(\begin{array}{c|ccc}
a_{11}^{-1} & -a_{11}^{-1}a_{12}&\cdots &-a_{11}^{-1}a_{1n} \\ \hline
 &1&&\\
&&\ddots&\\
 &&& 1
\end{array}\right)\in GL_n(R)
$$
such that
$$AS_1=\left(\begin{array}{cc} 1 & 0\\B_1& B_2\end{array}\right),$$
where $B_1=(a_{11}^{-1}a_{21},\ldots,a_{11}^{-1}a_{m1})^\top$ and $B_2\in R^{(m-1)\times(n-1)}$. By Proposition~\ref{prop1.1} (ii), ${\rm rk}(AS_1)=m$,
which implies that the set of all the rows of $AS_1$ is linearly independent.
It follows that $\{\beta_2,\ldots,\beta_m\}$ is linearly independent, where
$B_2=(\beta_2^{\top},\ldots,\beta_m^{\top})^{\top}$. By Proposition~\ref{prop1.1} (iii), ${\rm rk}(B_2)=m-1$.
By induction, there exists  an $S_2\in GL_{n-1}(R)$ such that
$$B_2S_2=\left(\begin{array}{cc} I_{m-1}& 0_{m-1,n-m}\end{array}\right).$$
Then $$AS_1\left(\begin{array}{cc} 1  & 0\\0& S_2\end{array}\right)
=\left(\begin{array}{ccc}1&0&0\\B_1& I_{m-1}& 0_{m-1,n-m}\end{array}\right).$$
Write $$
S_3=\left(\begin{array}{cc|c}
1&0&\\
-B_1& I_{m-1} &\\ \hline
&&I_{n-m}\end{array}\right).$$
Then
$$S=S_1\left(\begin{array}{cc}1&0\\0& S_2\end{array}\right)S_3\in GL_n(R)$$
and $AS=(I_m,0_{m,n-m})$.
Therefore, the desired result follows. $\qed$

\begin{lemma}\label{lemma2.4}
Let $m\leq n$ and $A=(\alpha_1^{\top},\alpha_2^{\top},\ldots,\alpha_m^{\top})^{\top}\in R^{m\times n}$. Then the following statements are equivalent.
\begin{itemize}
\item[\rm(i)]
${\rm rk}(A)=m$.

\item[\rm(ii)]
$\{\alpha_1,\alpha_2,\ldots,\alpha_m\}$ is  linearly independent.

\item[\rm(iii)]
There exists an $S\in GL_n(R)$ such that
$AS=(I_m,0_{m,n-m})$.

\item[\rm(iv)]
$\{\alpha_1,\alpha_2,\ldots,\alpha_m\}$ is unimodular.
\end{itemize}
\end{lemma}
\proof (i) $\Leftrightarrow$ (ii). This is a direct result of Proposition~\ref{prop1.1} (iii).

(i) $\Leftrightarrow$ (iii). This is a direct result of Lemma~\ref{lemma2.3}.

(iii) $\Leftrightarrow$ (iv). If there exists an $S\in GL_n(R)$ such that
$AS=(I_m,0_{m,n-m})$, then $A=(I_m,0_{m,n-m})S^{-1}$. Pick $$B=S\left(\begin{array}{c}I_m\\0_{n-m,m}\end{array}\right).$$ Then
$AB=I_m$. So, $\{\alpha_1,\alpha_2,\ldots,\alpha_m\}$ is unimodular.

If $\{\alpha_1,\alpha_2,\ldots,\alpha_m\}$ is unimodular, then there exists some $B\in R^{n\times m}$ such that $AB=I_m$.
By (\ref{shi1}), $\pi(I_m)=\pi(AB)=\pi(A)\pi(B)$, which implies that ${\rm rk}(\pi(A))=m$ since $R/M$ is a field.
By Lemma~\ref{lemma2.1}, ${\rm rk}(A)=m$. So, there exists an $S\in GL_n(R)$ such that
$AS=(I_m,0_{m,n-m})$ by Lemma~\ref{lemma2.3}. $\qed$

 For $T\in GL_m(R),S\in GL_n(R)$ and ${\cal M}\subseteq R^{m\times n}(m)$, we write
$T{\cal M}=\{TA: A\in {\cal M}\}$ and ${\cal M}S=\{AS: A\in {\cal M}\}.$

\begin{lemma}\label{lemma2.5}
Let $0\leq m\leq n$. Then
$$|R^{m\times n}(m)|=|R|^{\frac{m(m-1)}{2}}\prod_{i=0}^{m-1}(|R|^{n-i}-|M|^{n-i}).$$
In particular,
$$|GL_{n}(R)|=|R|^{\frac{n(n-1)}{2}}\prod_{i=0}^{n-1}(|R|^{n-i}-|M|^{n-i}).$$
\end{lemma}
\proof Let $A=(\alpha_1^{\top},\alpha_2^{\top},\ldots,\alpha_m^{\top})^{\top}\in R^{m\times n}(m)$. We enumerate how many $A$'s there are.
 By Lemma~\ref{lemma2.4},
$\{\alpha_1,\alpha_2,\ldots,\alpha_m\}$ is linearly independent, which implies that $\{\alpha_1\}$ is  linearly independent.
By Lemma~\ref{lemma2.2}, there are $|R|^n-|M|^n$ choices for $\alpha_1$. Once $\alpha_1$ is chosen,
let $$R^{m\times n}(m,\alpha_1)=\left\{\left(\begin{array}{c}\alpha_1\\ A_2\end{array}\right)
\in R^{m\times n}(m): A_2\in R^{(m-1)\times n}(m-1)\right\}.$$
Similar to the proof of Lemma~\ref{lemma2.3}, there exists
$S_1\in GL_n(R)$ such that $R^{m\times n}(m,\alpha_1)S_1=R^{m\times n}(m,e_1)$, where
$e_1$ is the row vector in $R^{n}$ in which $1$st coordinate is 1 and all other
coordinates are $0$. Then $|R^{m\times n}(m,\alpha_1)|=|R^{m\times n}(m,e_1)|$, which implies that
 $|R^{m\times n}(m,\alpha_1)|$ is independent of the particular choice
of $\alpha_1$. Without loss of generality, we may assume that
$\alpha_1=e_1$ and $A_2=(\alpha_2^\top,\ldots,\alpha_m^\top)^\top$. Since $\{e_1,\alpha_2\}$ is linearly independent, by Lemma~\ref{lemma2.2} again,
there are $|R|(|R|^{n-1}-|M|^{n-1})$ choices for $\alpha_2$. Proceeding in this way, finally after $\alpha_1,\ldots,\alpha_{m-1}$
are chosen, there are $|R|^{m-1}(|R|^{n-m+1}-|M|^{n-m+1})$ choices for $\alpha_m$. Therefore,
$$|R^{m\times n}(m)|=\prod_{i=0}^{m-1}|R|^{i}(|R|^{n-i}-|M|^{n-i})=|R|^{\frac{m(m-1)}{2}}\prod_{i=0}^{m-1}(|R|^{n-i}-|M|^{n-i}),$$
as desired. $\qed$

Recall that $R^{m\times n}(m,A_1)$ is as in (\ref{equa-zj}).
By the proof of Lemma~\ref{lemma2.5}, we obtain the following result.

 \begin{cor}\label{lemma2.6}
Let $0\leq m_1\leq m\leq n$ and $A_1\in R^{m_1\times n}(m_1)$ be a fixed matrix. Then
$$|R^{m\times n}(m,A_1)|=|R|^{\frac{(m-m_1)(m+m_1-1)}{2}}\prod_{i=m_1}^{m-1}(|R|^{n-i}-|M|^{n-i}).$$
 \end{cor}

\section{Matrices over commutative rings}
In this section, we always assume that $R$ is a commutative ring. It is well known that $R$ is a product
of local rings (cf. \cite{Atiyah}). Write
$$(R,+,\cdot)=(R_1\times R_2\times\cdots\times R_\ell,+,\cdot),$$
where $R_1,R_2,\ldots,R_\ell$ are local rings  with maximal ideals $M_1, M_2,\ldots, M_\ell$, respectively.
Let $$\begin{array}{rccl}
\rho_i: & R & \longrightarrow & R_i,\\
&(r_1,r_2,\ldots,r_\ell) & \longmapsto & r_i
\end{array}$$ be the projection map for all $i\in [\ell]=\{1,2,\ldots,\ell\}$.
For $A=(a_{uv})\in R^{m\times n}$ and $B\in R^{n\times k}$, let
$\rho_i(A)=(\rho_i(a_{uv}))$ for all $i\in[\ell]$.
Since $\rho_i(x+y)=\rho_i(x)+\rho_i(y)$ and
$\rho_i(xy)=\rho_i(x)\rho_i(y)$ for all $x,y\in R$, it is easy to prove that
 \begin{equation}\label{shi4}
 \rho_i(AB)=\rho_i(A)\rho_i(B)\;\hbox{for all}\;i\in[\ell].
 \end{equation}
Therefore, we have
$$\rho_i(R^{m\times n}):=\{\rho_i(A) : A\in R^{m\times n}\}=R_i^{m\times n},\;\rho_i(GL_n(R))=GL_n(R_i)\quad\hbox{for all}\;i\in[\ell].$$
From \cite{McDonald}, we deduce that
\begin{equation}\label{shi2}
(GL_n(R),\cdot)\cong (GL_{n}(R_1)\times GL_{n}(R_2)\times\cdots\times GL_n(R_\ell),\cdot)
\end{equation}
and \begin{equation}\label{shi3}
(R^{m\times n},+)\cong (R_1^{m\times n}\times R_2^{m\times n}\times \cdots\times R_\ell^{m\times n},+).
\end{equation}

\begin{lemma}\label{lemma3.0}{\rm (See \cite{Bollman}).}
If $A\in R^{m\times n}$, then ${\rm rk}(A)=\min\{{\rm rk}(\rho_i(A)) : 1\leq i\leq\ell\}$.
\end{lemma}

\begin{lemma}\label{lemma3.1}{\rm (See \cite[Proposition~1.3]{Sirisuk}).}
Let $x_1,x_2,\ldots,x_s\in R^n$. Then $\{x_1,x_2,\ldots,x_s\}$ is linearly independent
over $R$ if and only if $\{\rho_i(x_1),\rho_i(x_2),\ldots,\rho_i(x_s)\}$ is linearly independent over $R_i$ for all $i\in[\ell]$.
\end{lemma}

For convenience, we write the identity (resp. the zero) of $R_i$ as $1$ (resp.  $0$) for all $i\in[\ell]$,
and the identity and the zero of $R$ also as $1$ and $0$,\footnote{In fact, the identity and the zero of $R$ are $\underbrace{(1,\ldots,1)}_\ell$ and $\underbrace{(0,\ldots,0)}_{\ell}$, respectively.} respectively.

\begin{lemma}\label{lemma3.2}
Let $m\leq n$ and $A\in R^{m\times n}$. Then ${\rm rk}(A)=m$ if and only if there exists an $S\in GL_n(R)$ such that
$AS=(I_m,0_{m,n-m})$.
\end{lemma}
\proof If  there exists an $S\in GL_n(R)$ such that
$AS=(I_m,0_{m,n-m})$, then ${\rm rk}(A)={\rm rk}(I_m,0_{m,n-m})=m$ by Proposition~\ref{prop1.1}.

Suppose ${\rm rk}(A)=m$. Let $A=(a_{ij})=(\alpha_1^{\top},\alpha_2^{\top},\ldots,\alpha_m^{\top})^{\top}$. By Proposition~\ref{prop1.1} (iii), $\{\alpha_1,\alpha_2,\ldots,\alpha_m\}$
 is linearly independent over $R$. By Lemma~\ref{lemma3.1}, $\{\rho_i(\alpha_1),\rho_i(\alpha_2),\ldots,\rho_i(\alpha_m)\}$
 is linearly independent over $R_i$ for all $i\in[\ell]$. For every $i\in[\ell]$, by Lemma~\ref{lemma2.4}, there exists an
 $S_i\in GL_n(R_i)$ such that $\rho_i(A)S_i=(I_m,0_{m,n-m})$. From (\ref{shi2}), we deduce that there exists the  unique
 $S\in GL_n(R)$ such that $\rho_i(S)=S_i$ for all $i\in[\ell]$. By (\ref{shi4}), we obtain
 $$\rho_i(AS)=\rho_i(A)\rho_i(S)=(I_m,0_{m,n-m})=\rho_i(I_m,0_{m,n-m})\;\hbox{for all}\;i\in[\ell].$$
 By (\ref{shi3}), $AS=(I_m,0_{m,n-m})$. $\qed$

\begin{lemma}\label{lemma3.3}
Let $m\leq n$ and $A=(\alpha_1^{\top},\alpha_2^{\top},\ldots,\alpha_m^{\top})^{\top}\in R^{m\times n}$. Then the following statements are equivalent.
\begin{itemize}
\item[\rm(i)]
${\rm rk}(A)=m$.

\item[\rm(ii)]
$\{\alpha_1,\alpha_2,\ldots,\alpha_m\}$ is  linearly independent.

\item[\rm(iii)]
There exists an $S\in GL_n(R)$ such that
$AS=(I_m,0_{m,n-m})$.

\item[\rm(iv)]
$\{\alpha_1,\alpha_2,\ldots,\alpha_m\}$ is unimodular.
\end{itemize}
\end{lemma}
\proof (i) $\Leftrightarrow$ (ii). This is a direct result of Proposition~\ref{prop1.1} (iii).

(i) $\Leftrightarrow$ (iii). This is a direct result of Lemma~\ref{lemma3.2}.

(iii) $\Rightarrow$ (iv). The proof is similar to that of Lemma~\ref{lemma2.3}, and will be omitted.

(iv) $\Rightarrow$ (iii). If $\{\alpha_1,\alpha_2,\ldots,\alpha_m\}$ is unimodular, then there exists some $B\in R^{n\times m}$ such that $AB=I_m$.
By (\ref{shi4}), $I_m=\rho_i(I_m)=\rho_i(AB)=\rho_i(A)\rho_i(B)$ for all $i\in[\ell]$,
which implies that $\{\rho_i(\alpha_1),\rho_i(\alpha_2),$ $\ldots,\rho_i(\alpha_m)\}$ is unimodular for all $i\in[\ell]$.
By Lemma~\ref{lemma2.4}, $\{\rho_i(\alpha_1),\rho_i(\alpha_2),\ldots,$ $\rho_i(\alpha_m)\}$
is linearly independent over $R_i$ for all $\in[\ell]$, which implies that $\{\alpha_1,\alpha_2,\ldots,\alpha_m\}$
is linearly independent over $R$ by Lemma~\ref{lemma3.1}.
From (ii) $\Leftrightarrow$ (iii), we deduce that there exists an $S\in GL_n(R)$ such that
$AS=(I_m,0_{m,n-m})$. $\qed$

By Lemmas~\ref{lemma2.4},~\ref{lemma3.1} and~\ref{lemma3.3}, we obtain the following result.

\begin{cor}\label{lemma3.4}
Let $m\leq n$ and $A=(\alpha_1^{\top},\alpha_2^{\top},\ldots,\alpha_m^{\top})^{\top}\in R^{m\times n}$. Then the following hold.
\begin{itemize}
\item[\rm(i)]
${\rm rk}(A)=m$ if and only if ${\rm rk}(\rho_i(A))=m$ for all $i\in[\ell]$.

\item[\rm(ii)]
$\{\alpha_1,\alpha_2,\ldots,\alpha_m\}$ is  linearly independent over $R$ if and only if
$\{\rho_i(\alpha_1),\rho_i(\alpha_2),\ldots,\rho_i(\alpha_m)\}$ is  linearly independent over $R_i$ for all $i\in[\ell]$.

\item[\rm(iii)]
There exists an $S\in GL_n(R)$ such that
$AS=(I_m,0_{m,n-m})$ over $R$ if and only if there exists an $S_i\in GL_n(R_i)$ such that
$\rho_i(A)S_i=(I_m,0_{m,n-m})$ over $R_i$ for all $i\in[\ell]$.

\item[\rm(iv)]
$\{\alpha_1,\alpha_2,\ldots,\alpha_m\}$ is unimodular over $R$ if and only if
$\{\rho_i(\alpha_1),\rho_i(\alpha_2),\ldots,\rho_i(\alpha_m)\}$ is unimodular  over $R_i$ for all $i\in[\ell]$.
\end{itemize}
\end{cor}

\begin{lemma}\label{lemma3.5}
Let $m\leq n$ and $\alpha_1,\alpha_2,\ldots,\alpha_m\in R^{n}$.
If $\{\alpha_1,\alpha_2,\ldots,\alpha_m\}$ is linearly independent, then $\{\alpha_1,\alpha_2,\ldots,\alpha_m\}$
can be extended to a  unimodular basis of $R^{n}$.
\end{lemma}
\proof Since $\{\alpha_1,\alpha_2,\ldots,\alpha_m\}$ is linearly independent, by Lemma~\ref{lemma3.3},
there is an $S\in GL_n(R)$ such that $A=(I_m,0_{m,n-m})S^{-1}$, where $A=(\alpha_1^{\top}, \alpha_2^{\top}, \ldots, \alpha_m^{\top})^{\top}$.
Then
 $$
\widetilde{A}=\left(\begin{array}{c}
A\\ (0_{n-m,m},I_{n-m})S^{-1}
\end{array}\right)
=S^{-1}$$
is invertible. Write $\widetilde{A}=(\alpha_1^{\top},\ldots, \alpha_m^{\top},\alpha_{m+1}^{\top},\ldots,\alpha_n^{\top})^{\top}$.
Then $\{\alpha_1,\ldots, \alpha_m,\alpha_{m+1},\ldots,\alpha_n\}$ is a  unimodular basis of $R^{n}$. $\qed$

Recall that $R^{m\times n}(m,A_1)$ is as in (\ref{equa-zj}).
By (\ref{shi2}), (\ref{shi3}), Lemma~\ref{lemma2.5}, Corollaries~\ref{lemma2.6} and~\ref{lemma3.4},
we obtain the following results.

\begin{lemma}\label{lemma3.6}
Let $0\leq m\leq n$. Then
$$|R^{m\times n}(m)|=|R|^{\frac{m(m-1)}{2}}\prod_{j=1}^{\ell}\prod_{i=0}^{m-1}(|R_j|^{n-i}-|M_j|^{n-i}).$$
In particular,
$$|GL_{n}(R)|=|R|^{\frac{n(n-1)}{2}}\prod_{j=1}^{\ell}\prod_{i=0}^{n-1}(|R_j|^{n-i}-|M_j|^{n-i}).$$
\end{lemma}

 \begin{cor}\label{lemma3.7}
Let $0\leq m_1\leq m\leq n$ and $A_1\in R^{m_1\times n}(m_1)$ be a fixed matrix. Then
$$|R^{m\times n}(m,A_1)|=|R|^{\frac{(m-m_1)(m+m_1-1)}{2}}\prod_{j=1}^{\ell}\prod_{i=m_1}^{m-1}(|R_j|^{n-i}-|M_j|^{n-i}).$$
 \end{cor}

  \section{Vector spaces}
  In this section, we always assume that $R=R_1\times R_2\times\cdots\times R_\ell$ is a commutative ring,
where $R_1,R_2,\ldots,R_\ell$ are local rings  with maximal ideals $M_1, M_2,\ldots, M_\ell$, respectively.

Let $\{\alpha_1,\alpha_2,\ldots,\alpha_m\}\subseteq R^{n}$ be linearly independent,
and $\langle\alpha_1,\alpha_2,\ldots,\alpha_m\rangle$ denote the $R$-module generated by $\alpha_1,\alpha_2,\ldots,\alpha_m$.
Then $X=\langle\alpha_1,\alpha_2,\ldots,\alpha_m\rangle$  is an $m$-subspace of $R^{n}$, and
the matrix $(\alpha_1^{\top}, \alpha_2^{\top}, \ldots, \alpha_m^{\top})^{\top}$
is called a {\it matrix representation} of $X$. The matrix representation is not unique. We use an $m\times n$ matrix $A$ with
${\rm rk}(A)=m$ to represent an $m$-subspace of $R^{n}$. For $A,B\in R^{m\times n}$ with ${\rm rk}(A)={\rm rk}(B)=m$,
both $A$ and $B$ represent the same $m$-subspace if and only if there exists a $T\in GL_m(R)$ such that $B=TA$.
For convenience, if $A$ is a matrix representation of an $m$-subspace $X$, then we write $X=A$.

For an $m$-subspace $X$ of $R^{n}$,  by Lemma~\ref{lemma3.3}, there exists an $S\in GL_n(R)$
such that $X$ has matrix representations
\begin{equation}\label{shi5}
X=(I_m,0_{m,n-m})S=T(I_m,0_{m,n-m})S\; \hbox{for all}\; T\in GL_m(R).
\end{equation}

There is an action of
$GL_{n}(R)$ on $R^{n}$ defined as
follows:
$$\begin{array}{ccc}
R^{n}\times GL_{n}(R) & \longrightarrow & R^{n},\\
((r_1,\ldots,r_n),S)&\longmapsto&(r_1,\ldots,r_{n})S.
\end{array}
$$
The above action induces an action on the set of subspaces of
$R^{n}$; i.e., a subspace $P$ is carried by $S\in
GL_{n}(R)$ to the subspace $PS$.

For $0\leq m\leq n$, let  ${\cal M}_R(m,n)$ denote the set of all the $m$-subspaces of $R^n$, and let
$N_R(m,n)$  denote the size of ${\cal M}_R(m,n)$. Sirisuk and Meemark \cite{Sirisuk} determined the size of ${\cal M}_R(m,n)$
by using  a lifting technique.

\begin{thm}\label{lem4.1}
${\cal M}_R(m,n)$ is non-empty if and only if $0\leq m\leq n$.
Moreover, if ${\cal M}_R(m,n)$ is non-empty, then it forms an orbit of subspaces under
the action of $GL_{n}(R)$ and
$$
N_R(m,n)=\prod_{j=1}^{\ell}\prod_{i=0}^{m-1}\frac{|R_j|^{n-i}-|M_j|^{n-i}}
{|R_j|^{m-i}-|M_j|^{m-i}}.
$$
\end{thm}
\proof Clearly, ${\cal M}_R(m,n)$ is non-empty if and only if $0\leq m\leq n$. For any $m$-subspace $P\in{\cal M}_R(m,n)$,
by Lemma~\ref{lemma3.3}, there exists an $S\in GL_n(R)$ such that $PS=(I_m,0_{m,n-m})$.
 By (\ref{shi5}),  ${\cal M}_R(m,n)$ forms an orbit of subspaces under the action of
$GL_{n}(R)$.

Recall that $R^{m\times n}(m)$ is the set of all the $m\times n$ matrices of McCoy rank $m$ over $R$.
By Lemma~\ref{lemma3.6} and (\ref{shi5}), $$N_R(m,n)=\frac{|R^{m\times n}(m)|}{|GL_m(R)|}=\prod_{j=1}^{\ell}\prod_{i=0}^{m-1}\frac{|R_j|^{n-i}-|M_j|^{n-i}}
{|R_j|^{m-i}-|M_j|^{m-i}},$$
as desired.
$\qed$

For a fixed $m$-subspace $P$ of $R^{n}$,
let ${\cal M}_R(m_1,m,n)$ denote the set of all the
$m_1$-subspaces contained in $P$, and let $ N_R(m_1,m,n)$
denote the size of  ${\cal M}_R(m_1,m,n)$.
  By Theorem~\ref{lem4.1}, $N_R(m_1,m,n)$ is independent of the particular choice
of the subspace $P$.

\begin{thm}\label{lem4.2}
${\cal M}_R(m_1,m,n)$ is non-empty if and only if $0\leq m_1\leq m\leq n$.
Moreover, if ${\cal M}_R(m_1,m,n)$ is non-empty, then
$$
N_R(m_1,m,n)=N_R(m_1,m)=\prod_{j=1}^{\ell}\prod_{i=0}^{m_1-1}\frac{|R_j|^{m-i}-|M_j|^{m-i}}
{|R_j|^{m_1-i}-|M_j|^{m_1-i}}.
$$
\end{thm}
\proof Suppose that $P$ is a fixed $m$-subspace of $R^{n}$. Then we may assume that $P=(I_m,0_{m,n-m})$.
Define $$\begin{array}{rccc}
\sigma:& P &\longrightarrow& R^m,\\
&(r_1,\ldots,r_m,\underbrace{0,\ldots,0}_{n-m}) &\longmapsto& (r_1,\ldots,r_m).
\end{array}$$
Then $\sigma$ is an isomorphism from $P$ to $R^m$. Therefore,
$N_R(m_1,m,n)=N_R(m_1,m)$. $\qed$

For a fixed $m_1$-subspace $P$ of $R^n$,
let ${\cal M}'_R(m_1,m,n)$ denote the set of all the
$m$-subspaces containing $P$, and let $N'_R(m_1,m,n)$ denote the size of ${\cal M}'_R(m_1,m,n)$.
  By Theorem~\ref{lem4.1}, $N'_R(m_1,m,n)$ is independent of the particular choice
of the subspace $P$.

\begin{thm}\label{lem4.3}
${\cal M}'_R(m_1,m,n)$ is non-empty if and only if $0\leq m_1\leq m\leq n$.
Moreover, if ${\cal M}'_R(m_1,m,n)$ is non-empty, then
$$N'_R(m_1,m,n)=N_R(m-m_1,n-m_1)=\prod_{j=1}^{\ell}\prod_{i=0}^{m-m_1-1}\frac{|R_j|^{n-m_1-i}-|M_j|^{n-m_1-i}}
{|R_j|^{m-m_1-i}-|M_j|^{m-m_1-i}}.$$
\end{thm}
\proof
By counting the number of couples $(P,Q)$, where $P\in{\cal M}_R(m_1,n)$ and $Q\in{\cal M}_R(m,n)$
  with $P\subseteq Q$, by Theorems~\ref{lem4.1} and~\ref{lem4.2},
  $$N'_R(m_1,m,n)=\frac{N_R(m,n)N_R(m_1,m,n)}{N_R(m_1,n)}=N_R(m-m_1,n-m_1).$$
Therefore, the desired result follows. $\qed$

Let $X$ and $Y$ be two subspaces of $R^{n}$.
The set of vectors belonging to both $X$ and $Y$ is a linear subset of $R^{n}$,
called the {\it intersection} of $X$ and $Y$ and denoted by $X\cap Y$.
The set of vectors which can be written as sums of a vector of $X$ and a vector of $Y$ is also
a linear subset of $R^{n}$, called the {\it join} of $X$ and $Y$ and denoted by $X+Y$.
In general, the intersection and the join of two subspaces in $R^{n}$ are not subspaces.
Note that $X\cap Y=Y\cap X,X+Y=Y+X$, and $X\cap Y=X,X+Y=Y$ if $X\subseteq Y$.
Clearly, $X\cap Y$ contains a subspace of dimension $\dim(X\cap Y)$.

\begin{thm}\label{lem4.4}{\rm(Dimensional formula).}
Let $A$ and $B$ be two subspaces of $R^n$. Then
$$
\max\{\dim(A),\dim(B)\}\leq\dim(A+B)\leq\min\{n,\dim(A) + \dim(B) - \dim(A\cap B)\},
$$
and the following statements are equivalent:
 \begin{itemize}
 \item[\rm(i)]
 $\dim(A+B)=\dim(A) + \dim(B) - \dim(A\cap B)$;

 \item[\rm(ii)]
  $A+B$ is a subspace of $R^n$;

 \item[\rm(iii)]
 $A\cap B$ is a subspace of $R^n$.
 \end{itemize}
\end{thm}
\proof Suppose $\dim(A\cap B)=t,\dim(A)=m_1,\dim(B)=m_2$ and $\dim(A+B)=k$. Then $k\leq n$.
Since $\dim(A\cap B)=t$, there is a largest unimodular vector subset $\{\alpha_1,\ldots,\alpha_t\}$
of $A\cap B$.  By Lemma~\ref{lemma3.5},
$\{\alpha_1,\ldots,\alpha_t\}$ can be extended to a  unimodular basis $\{\alpha_1,\ldots,\alpha_{m_1}\}$ of $A$ and
a  unimodular basis $\{\alpha_1,\ldots,\alpha_t,\beta_{t+1},\ldots,\beta_{m_2}\}$ of $B$.
Then
\begin{equation}\label{shi6}
A+B=\langle\alpha_1,\ldots,\alpha_{m_1},\beta_{t+1},\ldots,\beta_{m_2}\rangle
\end{equation} is a linear subset.
Let $\{\gamma_1,\ldots,\gamma_k\}$ be a largest unimodular vector subset
of $A+B$. Write $A_1=(\alpha_1^{\top},\ldots,\alpha_{m_1}^{\top},\beta_{t+1}^{\top},\ldots,\beta_{m_2}^{\top})^{\top}$
and $A_2=(\gamma_1^{\top},\ldots,\gamma_k^{\top})^{\top}$.
Since $\{\gamma_1,\ldots,\gamma_k\}\subseteq A+B$, there exists a $k\times(m_1+m_2-t)$ matrix $M$
such that $A_2=MA_1$. Since $\{\gamma_1,\ldots,\gamma_k\}$ is unimodular, by Lemma~\ref{lemma3.3}, there is an $n\times k$ matrix $D_2$ such that $A_2D_2=I_k$. Therefore, we have $M(A_1D_2)=I_k$,  which implies that ${\rm rk}(M)=k$ by Lemma~\ref{lemma3.3}.
By Proposition~\ref{prop1.1}, $k\leq m_1+m_2-t$. So, $k\leq\min\{n,m_1+m_2-t\}$.
Since $A\subseteq A+B$ and $B\subseteq A+B$, $\dim(A+B)\geq\max\{\dim(A),\dim(B)\}$.

(i) $\Leftrightarrow$ (ii).
Suppose $k=m_1+m_2-t$. Then $M$ is invertible, which implies that $M^{-1}A_2=A_1$.
If $XA_1=0_{1,n}$ where $X=(r_1,r_2,\ldots,r_{k})$, then $XM^{-1}A_2=0_{1,n}$,
which implies that $XM^{-1}=0_{1,k}$ since $A_2D_2=I_k$. Therefore, $X=0_{1,k}$ and
$\{\alpha_1,\ldots,\alpha_{m_1},\beta_{t+1},\ldots,\beta_{m_2}\}$ is linearly independent.
 From (\ref{shi6}), we deduce that $A+B$ is a $k$-subspace.

 Conversely, suppose that $A+B$ is a $k$-subspace of $R^n$. By Theorem~\ref{lem4.1},
 there is an $S\in GL_n(R)$ such that $(A+B)S=(I_k,0_{k,n-k})$. So, we may assume that
 $$AS=\left(\begin{array}{c}
 \alpha_1\\
 \vdots\\
 \alpha_{m_1}
 \end{array}\right)S=(I_{m_1},0_{m_1,n-m_1}),$$
 and $$BS=\left(\begin{array}{c}
 \alpha_1\\
 \vdots\\
 \alpha_t\\
 \beta_{t+1}\\
 \vdots\\
 \beta_{m_2}
 \end{array}\right)S=\left(\begin{array}{ccc}
 I_t & 0_{t,m_1-t} & 0_{t,n-m_1}\\
 0_{m_2-t,t} & B_2 & B_3
 \end{array}\right),$$
 where $B_2\in R^{(m_2-t)\times(m_1-t)},B_3\in R^{(m_2-t)\times(n-m_1)}$ and ${\rm rk}(B_2,B_3)=m_2-t$.
  Then
 $$
 \left(\begin{array}{ccc}
 I_{t} &  & \\
      & I_{m_1-t} & \\
  & -B_2 & I_{m_2-t}
 \end{array}\right)A_1S=\left(\begin{array}{cc}
 I_{m_1}  & 0_{m_1,n-m_1}\\
 0_{m_2-t,m_1}  & B_3
 \end{array}\right).
$$
 Write $B_3=(B_{31},B_{32})$, where $B_{31}\in R^{(m_2-t)\times(k-m_1)}$ and $B_{32}\in R^{(m_2-t)\times(n-k)}$.
 Since all the rows of $(0_{m_2-t,m_1},B_3)$ are contained in the subspace $(I_k,0_{k,n-k})$,
 there exists an $(m_2-t)\times k$ matrix $M_1$
such that $(0_{m_2-t,m_1},B_3)=M_1(I_k,0_{k,n-k})$, which implies that $B_{32}=0_{m_2-t,n-k}$.
Similarly,
 there exist a $(k-m_1)\times m_1$ matrix $M_2$ and a $(k-m_1)\times(m_2-t)$ matrix $M_3$
such that $$(0_{k-m_1,m_1},I_{k-m_1},0_{k-m_1,n-k})
=(M_2,M_3)\left(\begin{array}{ccc}
 I_{m_1}  & 0_{m_1,k-m_1}  & 0_{m_1,n-k}\\
 0_{m_2-t,m_1}  & B_{31} & 0_{m_2-t,n-k}
 \end{array}\right),$$ which implies that $M_2=0_{k-m_1,m_1}$ and $M_3B_{31}=I_{k-m_1}$.
 By Lemma~\ref{lemma3.3}, there exists a  $T\in GL_{m_2-t}(R)$ such that
 $$TB_{31}=\left(\begin{array}{c}
 I_{k-m_1}\\ 0_{m_1+m_2-t-k,k-m_1}\end{array}\right).
 $$
 Then
 $$
 \left(\begin{array}{cc}
 I_t \\
   & T
 \end{array}\right)BS
 =\left(\begin{array}{cccc}
 I_t & 0_{t,m_1-t} & 0_{t,k-m_1} & 0_{t,n-k}\\
 0_{m_2-t,t} & B_{21} & I_{k-m_1} & 0_{k-m_1,n-k}\\
 0           &  B_{22}     & 0_{m_1+m_2-t-k,k-m_1} & 0
 \end{array}\right),
 $$
 where $TB_2=\left(\begin{array}{c}B_{21}\\ B_{22}\end{array}\right)$. From ${\rm rk}(B_2,B_3)=m_2-t$, we deduce that
 ${\rm rk}(B_{22})=m_1+m_2-t-k$.
It follows that $t=\dim(A\cap B)=\dim(AS\cap BS)\geq m_1+m_2-k$. Since $k\leq m_1+m_2-t$, $k=m_1+m_2-t$.

 (i) $\Leftrightarrow$ (iii).
Suppose $k=m_1+m_2-t$. From (\ref{shi6}), we deduce that $\{\alpha_1,\ldots,\alpha_{m_1},\beta_{t+1},\ldots,\beta_{m_2}\}$
is a unimodular basis of $A+B$.
For any $\alpha\in A\cap B$, there exist $r_1,\ldots,r_{m_1},s_1,\ldots,s_{m_2}\in R$ such that
  $$\alpha=\sum_{j=1}^{m_1}r_j\alpha_j=\sum_{j=1}^{t}s_j\alpha_j+\sum_{j=t+1}^{m_2}s_j\beta_j.$$
It follows that  $\sum_{j=1}^{t}(r_j-s_j)\alpha_j+\sum_{j=t+1}^{m_1}r_j\alpha_j-\sum_{j=t+1}^{m_2}s_j\beta_j=0_{1,n}.$
Since $\{\alpha_1,\ldots,\alpha_{m_1},\beta_{t+1},$ $\ldots,\beta_{m_2}\}$ is linearly independent, $r_j=s_j$ for all $1\leq j\leq t$, $r_j=0$ for all $t+1\leq j\leq m_1$ and $s_j=0$ for all $t+1\leq j\leq m_2$.
So, $\alpha=\sum_{j=1}^{t}r_j\alpha_j\in\langle\alpha_1,\ldots,\alpha_{t}\rangle$.
Therefore, $A\cap B=\langle\alpha_1,\ldots,\alpha_{t}\rangle$ is a $t$-subspace of $R^n$.

Conversely, suppose that $A\cap B$ is a $t$-subspace of $R^n$. Then $\{\alpha_1,\ldots,\alpha_t\}$ is a unimodular basis of $A\cap B$.
If $\sum_{i=1}^{m_1}r_i\alpha_i+\sum_{j=t+1}^{m_2}s_j\beta_j=0_{1,n}$, then
$\sum_{i=1}^{m_1}r_i\alpha_i=-\sum_{j=t+1}^{m_2}s_j\beta_j\in A\cap B$, and there exist $s_1,\ldots,s_t\in R$ such that
$\sum_{i=1}^{t}s_i\alpha_i+\sum_{j=t+1}^{m_2}s_j\beta_j=0_{1,n}$. Since $\{\alpha_1,\ldots,\alpha_t,\beta_{t+1},\ldots,\beta_{m_2}\}$
is a  unimodular basis of $B$, $s_1=\cdots=s_{m_2}=0$. It follows that $\sum_{i=1}^{m_1}r_i\alpha_i=0_{1,n}$, which implies that
$r_1=\cdots=r_{m_1}=0$ since $\{\alpha_1,\ldots,\alpha_{m_1}\}$ is a  unimodular basis of $A$.
Therefore, $\{\alpha_1,\ldots,\alpha_{m_1},\beta_{t+1},\ldots,\beta_{m_2}\}$ is linearly independent.
By (\ref{shi5}), $A+B$ is an $(m_1+m_2-t)$-subspace and $k=m_1+m_2-t$.  $\qed$

\begin{cor}\label{lem4.5}
Let $A$ and $B$ be two subspaces of $R^n$. Then $A\cap B=\{0_{1,n}\}$ if and only if $A+B$ is a $(\dim(A)+\dim(B))$-subspace of $R^n$,
and $A+B=R^n$ if and only if $A\cap B$ is a $(\dim(A)+\dim(B)-n)$-subspace of $R^n$.
\end{cor}

\medskip
\noindent{\bf Remarks.}
Let $A$ and $B$ be two subspaces of $R^n$. If $R$ is a field, then
$\dim(A+B)=\dim(A) + \dim(B) - \dim(A\cap B)$ always holds. Suppose that $R$ is not a field.
Since $R$ is a finite ring, each non-zero element in $R\setminus R^\ast$ is a zero divisor.
Let $1\leq m_1\leq m_2, m_1+m_2\leq n$ and $a\in R\setminus R^\ast$ be a fixed non-zero element.
Pick $A=(0_{m_1,n-m_1-m_2},aI_{m_1},0_{m_1,m_2-m_1},I_{m_1})$ and $B=(0_{m_2,n-m_2},I_{m_2})$.
Then $\dim(A+B)=m_2$ and $\dim(A\cap B)=0$. So, $\dim(B)=m_2=\dim(A+B)<m_1+m_2=\dim(A)+\dim(B)-\dim(A\cap B)$.

\medskip
Let $A$ be an $m$-subspace of $R^n$, and
$$
A^\perp=\{y\in R^n: yx^{\top}=0\;\hbox{for all}\;x\in A\}.
$$
Then $A^\perp$ is a linear subset of $R^n$.  By Lemma~\ref{lemma3.3}, $A$ has a matrix representation
$A=(I_m,0_{m,n-m})S$, where $S\in GL_n(R)$. It is easy to prove that $A^\perp$  has a matrix
representation
\begin{equation}\label{shi-N5}
A^\perp=(0_{n-m,m},I_{n-m})(S^{-1})^{\top}\; \hbox{if}\; A=(I_m,0_{m,n-m})S.
\end{equation}
Hence $A^\perp$ is an $(n-m)$-subspace of $R^n$.
The subspace $A^\perp$ is called the {\it dual subspace} of $A$. Note that
$\dim(A)+\dim(A^\perp)=n$ and $(A^{\perp})^\perp=A.$
If $A_1$ and $A_2$ are two subspaces of $R^n$, then $A_1\subseteq A_2$ if and only if $A_2^\perp\subseteq A_1^\perp.$

\begin{thm}\label{lem4.6}
Let $A$ and $B$ be two subspaces of $R^n$.  If $\dim(A+B)=\dim(A) + \dim(B) - \dim(A\cap B)$, then
$$(A\cap B)^{\bot}=A^{\bot} +B^{\bot}\quad\hbox{and}\quad(A+B)^{\bot}=A^{\bot}\cap B^{\bot}.$$
\end{thm}
\proof Since $\dim(A+B)=\dim(A) + \dim(B) - \dim(A\cap B)$, by Theorem~\ref{lem4.4}, both $A\cap B$ and $A+B$ are subspaces of $R^n$.
By (\ref{shi-N5}), $A^{\bot},B^{\bot},(A\cap B)^{\bot}$ and $(A+B)^{\bot}$ are subspaces.
Suppose $\dim(A)=m_1,\dim(B)=m_2$ and $\dim(A\cap B)=t$. Then
$\dim(A+B)=m_1+m_2-t,\dim(A^{\bot})=n-m_1,\dim(B^{\bot})=n-m_2,\dim((A\cap B)^{\bot})=n-t$
and  $\dim((A+B)^{\bot})=n+t-m_1-m_2$. Write $A=((A\cap B)^{\top},A_1^{\top})^{\top}$
and $B=((A\cap B)^{\top},B_1^{\top})^{\top}$.

By Theorems~\ref{lem4.1} and~\ref{lem4.4}, there is an $S\in GL_n(R)$ such that
$$(A+B)S^{-1}=
\left(\begin{array}{c}
A\cap B \\
A_1\\
B_1
\end{array}\right)S^{-1}=\left(\begin{array}{cccc}
I_{t}& 0 & 0 & 0_{t,n+t-m_1-m_2}\\
0 & I_{m_1-t} & 0 & 0\\
0&0& I_{m_2-t} & 0
\end{array}\right).$$ It follows that $A=(I_{m_1},0_{m_1,n-m_1})S$
and $$B=\left(\begin{array}{cccc}
I_{t}& 0_{t,m_1-t} & 0 & 0_{t,n+t-m_1-m_2}\\
0&0& I_{m_2-t} & 0
\end{array}\right)S.$$
By (\ref{shi-N5}) again, $A^{\bot}=(0_{n-m_1,m_1},I_{n-m_1})(S^{-1})^{\top}$
and
$$B^{\bot}=\left(\begin{array}{cccc}
0_{m_1-t,t}& I_{m_1-t} & 0 & 0\\
0&0& 0_{n+t-m_1-m_2,m_2-t} & I_{n+t-m_1-m_2}
\end{array}\right)(S^{-1})^{\top}.$$
Therefore
$$A^{\bot}+B^{\bot}=\left(\begin{array}{ccc}
0_{m_1-t,t}& I_{m_1-t} &  0\\
0&0 & I_{n-m_1}
\end{array}\right)(S^{-1})^{\top}$$
is an $(n-t)$-subspace and $\dim(A^{\bot}\cap B^{\bot})\geq n+t-m_1-m_2$.

From $A^{\bot},B^{\bot}\subseteq(A\cap B)^{\bot}$, we deduce that $A^{\bot}+B^{\bot}\subseteq(A\cap B)^{\bot}$.
Since both $A^{\bot}+B^{\bot}$ and $(A\cap B)^{\bot}$ are $(n-t)$-subspaces,
$(A\cap B)^{\bot}=A^{\bot} +B^{\bot}$.
By Theorem~\ref{lem4.4} and $\dim(A^{\bot}\cap B^{\bot})\geq n+t-m_1-m_2$,
$\dim(A^{\bot}\cap B^{\bot})=n+t-m_1-m_2$, and therefore
$$A^{\bot}\cap B^{\bot}=(0_{n+t-m_1-m_2,m_1+m_2-t},I_{n+t-m_1-m_2})(S^{-1})^{\top}$$ is an $(n+t-m_1-m_2)$-subspace.
From $(A+B)^{\bot}\subseteq A^{\bot},B^{\bot}$, we deduce that $(A+B)^{\bot}\subseteq A^{\bot}\cap B^{\bot}$.
Since both $(A+B)^{\bot}$ and $A^{\bot}\cap B^{\bot}$ are $(n+t-m_1-m_2)$-subspaces, $(A+B)^{\bot}=A^{\bot}\cap B^{\bot}.$
Hence, the desired results follow. $\qed$

\section{Singular linear spaces}
 In this section, we always assume that $R=R_1\times R_2\times\cdots\times R_\ell$ is a commutative ring,
where $R_1,R_2,\ldots,R_\ell$ are local rings  with maximal ideals $M_1, M_2,\ldots, M_\ell$, respectively.

For two non-negative integers $n$ and $k$,
$R^{n+k}$   denotes the $(n+k)$-dimensional   row vector
space  over $R$. The set of all $(n+k)\times(n+k)$
invertible matrices over $R$ of the form
$$
\left(\begin{array}{cc}
 T_{11} & T_{12}\\
0 & T_{22}
\end{array}\right),$$
where $T_{11}\in GL_n(R)$ and $T_{22}\in GL_k(R)$, forms a group under matrix
multiplication, called the {\it singular general linear group} of
degree $n+k$ over $R$ and denoted by $GL_{n+k,n}(R)$.
There is an action of
$GL_{n+k,n}(R)$ on $R^{n+k}$ defined as
follows:
$$
\begin{array}{ccc}
R^{n+k}\times GL_{n+k,n}(R) &\longrightarrow& R^{n+k},\\
((x_1,x_2,\ldots,x_{n+k}),T)&\longmapsto&
(x_1,x_2,\ldots,x_{n+k})T.
\end{array}
$$
The above action induces an action on the set of subspaces of
$R^{n+k}$; i.e., a subspace $P$ is carried by $T\in
GL_{n+k,n}(R)$ to the subspace $PT$. The vector space
$R^{n+k}$ together with the above group action, is called
the $(n+k)$-dimensional {\it singular linear space} over
$R$.

For $1\leq i\leq n+k$, let $e_i$ denote the row vector in
$R^{n+k}$ in which $i$th coordinate is 1 and all other
coordinates are 0. Define $E=\langle e_{n+1},e_{n+2},\ldots,e_{n+k}\rangle$.
Then $E$ is a $k$-subspace of $R^{n+k}$.
An $m$-subspace $P$ of $R^{n+k}$ is called an
 $(m,t)$-subspace if $P\cap E$ is a $t$-subspace. The collection
of all the $(m,0)$-subspaces of $R^{n+k}$, where
$0\leq m\leq n,$ is the {\em attenuated space}.

Let ${\cal M}_R(m,t;n+k,n)$ denote the set of all the $(m,t)$-subspaces of
 $R^{n+k}$, and let  $N_R(m,t;n+k,n)$
denote the size of ${\cal M}_R(m,t;n+k,n)$.

\begin{thm}\label{lem5.1}
${\cal M}_R(m,t;n+k,n)$ is non-empty if and only if
$0\leq t\leq k$ and $0\leq m-t\leq n$. Moreover, if ${\cal M}_R(m,t;n+k,n)$ is non-empty,
then it forms an orbit of subspaces under the action of $GL_{n+k,n}(R)$ and
$$
N_R(m,t;n+k,n)=|R|^{(m-t)(k-t)}N_R(m-t,n)N_R(t,k).
$$
\end{thm}
\proof If $0\leq t\leq k$ and $0\leq m-t\leq n$, then
$\langle e_1,\ldots,e_{m-t},e_{n+1},\ldots,e_{n+t}\rangle$ is an $(m,t)$-subspace, which implies that
${\cal M}_R(m,t;n+k,n)$ is non-empty. Conversely, suppose $P\in {\cal M}_R(m,t;n+k,n)$.
Since $P\cap E$ is a $t$-subspace of $E$, it has
a matrix representation $(0_{t,n},P_{22})$, where $P_{22}\in R^{t\times k}(t)$.
So, $0\leq t\leq k$.
By Lemma~\ref{lemma3.5}, $P$ has a matrix representation
\begin{equation}\label{shi7}
P=\left(\begin{array}{cc}
P_{11}&P_{12}\\
0_{t,n}&P_{22}
\end{array}\right),
\end{equation}
 where $P_{11}\in R^{(m-t)\times n}$ and $P_{12}\in R^{(m-t)\times k}$.
 By Theorem~\ref{lem4.4} and (\ref{shi5}), $P+E$ is an $(m-t+k)$-subspace  and has a matrix representation
 $$\left(\begin{array}{cc}
P_{11}&0_{m-t,k}\\
0_{t,n}&I_k
\end{array}\right),
 $$
 where $P_{11}\in R^{(m-t)\times n}(m-t)$.
 So, $0\leq m-t\leq n$.

For each $P\in {\cal M}_R(m,t;n+k,n)$  has a matrix representation as in (\ref{shi7}) with $P_{11}\in R^{(m-t)\times n}(m-t)$
and $P_{22}\in R^{t\times k}(t)$.
By Lemma~\ref{lemma3.3} and (\ref{shi5}), $P_{11}$ is an $(m-t)$-subspace
of $R^{n}$, and $P_{22}$ is a $t$-subspace
of $E$. By Theorem~\ref{lem4.1}, there exist  $T_{11}\in
GL_n(R)$ and $T_{22}\in GL_k(R)$ such that $P_{11}T_{11}=(I_{m-t},0_{m-t,n-m+t})$ and
$P_{22}T_{22}=(I_t,0_{t,k-t})$. Let
$$
T_1=\left(\begin{array}{cc}
 T_{11} & 0\\
 0& T_{22}
\end{array}\right).
$$
Then $PT_1$ has a matrix representation
$$PT_1=\left(\begin{array}{cccc}
 I_{m-t} & 0_{m-t,n-m+t} & 0 & P_{14}\\
 0& 0 & I_t & 0_{t,k-t}
\end{array}\right).$$
Let
$$
T_2=\left(\begin{array}{ccc|c}
 I_{m-t} & & & -P_{14}\\ \hline
 & I_{n-m+t} &&\\
 && I_t&\\
 &&&I_{k-t}
\end{array}\right).
$$
Then $T_1T_2\in GL_{n+k,n}(R)$ and
$$
PT_1T_2=\left(\begin{array}{cccc}
 I_{m-t} & 0_{m-t,n-m+t} & 0 & 0\\
 0& 0 & I_t & 0_{t,k-t}
\end{array}\right).
$$
Therefore,
${\cal M}_R(m,t;n+k,n)$   forms an orbit of subspaces under the action of
$GL_{n+k,n}(R)$.

Denote by $n_R(m,t;n+k,n)$ the number of all the matrices $P$ of the form (\ref{shi7})
with  $P_{11}\in R^{(m-t)\times n}(m-t)$
and $P_{22}\in R^{t\times k}(t)$.
Suppose that $P$ and $Q$ represent the same subspace, then there is
an $m\times m$ invertible matrix $S$ such that $SP=Q$. It follows
that $S$ is necessarily of the form
$$
S=\left(\begin{array}{cc}
 S_{11} & S_{12}\\
 0& S_{22}
\end{array}\right)\in GL_{m,m-t}(R).
$$
Moreover, if $SP=P$, then $S=I_m$. Consequently,
$$
N_R(m,t;n+k,n)=\frac{n_R(m,t;n+k,n)}{|GL_{m,m-t}(R)|}.
$$
By Lemma~\ref{lemma3.6} and Theorem~\ref{lem4.1}, the desired result follows.
$\qed$

For a fixed $(m,t)$-subspace $P$ of $R^{n+k}$,
let ${\cal M}_R(m_1,t_1;m,t;n+k,n)$ denote the set of all the
$(m_1,t_1)$-subspaces contained in $P$, and let $N_R(m_1,t_1;m,t;n+k,n)$
denote the size of  ${\cal M}_R(m_1,t_1;m,t;n+k,n)$.
  By Theorem~\ref{lem5.1}, $N_R(m_1,t_1;m,t;n+k,n)$ is independent of the particular choice
of the subspace $P$.

\begin{thm}\label{lem5.2}
${\cal M}_R(m_1,t_1;m,t;n+k,n)$ is non-empty if and only if
\begin{equation}\label{shi8}
0\leq t_1\leq t\leq k\quad and\quad 0\leq m_1-t_1\leq m-t\leq n.
\end{equation}
Moreover, if ${\cal M}_R(m_1,t_1;m,t;n+k,n)$ is non-empty. Then
$$
N_R(m_1,t_1;m,t;n+k,n)=|R|^{(m_1-t_1)(t-t_1)}N_R(m_1-t_1,m-t)N_R(t_1,t).
$$
\end{thm}
\proof Suppose that $P$ is a fixed $(m,t)$-subspace of $R^{n+k}$. Then we may assume that
$$P=\left(\begin{array}{cccc}
 I_{m-t} & 0_{m-t,n-m+t} & 0 & 0\\
 0& 0 & I_t & 0_{t,k-t}
\end{array}\right).$$
Define $$\begin{array}{rccc}
\sigma:& P &\longrightarrow& R^m,\\
&(r_1,\ldots,r_{m-t},\underbrace{0,\ldots,0}_{n-m+t},r_{n+1},\ldots,r_{n+t},\underbrace{0,\ldots,0}_{k-t}) &\longmapsto& (r_1,\ldots,r_{m-t},r_{n+1},\ldots,r_{n+t}).
\end{array}$$
Then $\sigma$ is an isomorphism from $P$ to $R^m$. So,
$${\cal M}_R(m_1,t_1;m,t;n+k,n)\not=\emptyset\Leftrightarrow {\cal M}_R(m_1,t_1;m,m-t)\not=\emptyset$$ and
$$N_R(m_1,t_1;m,t;n+k,n)=N_R(m_1,t_1;m,m-t).$$ By Theorem~\ref{lem5.1}, the desired result follows. $\qed$

For a fixed $(m_1,t_1)$-subspace $P$ of $R^{n+k}$,
let ${\cal M}'_R(m_1,t_1;m,t;n+k,n)$ denote the set of all the
$(m,t)$-subspaces containing $P$, and let $N'_R(m_1,t_1;m,t;n+k,n)$ denote the size of ${\cal M}'_R(m_1,t_1;m,t;n+k,n)$.
  By Theorem~\ref{lem5.1},  $N'_R(m_1,t_1;m,t;n+k,n)$ is independent of the particular choice
of the subspace $P$.

\begin{thm}\label{lem5.3}
${\cal M}'_R(m_1,t_1;m,t;n+k,n)$ is non-empty if and only if
(\ref{shi8}) holds.
Moreover, if ${\cal M}'_R(m_1,t_1;m,t;n+k,n)$ is non-empty, then
$$N'_R(m_1,t_1;m,t;n+k,n)=|R|^{(k-t)(m-t-m_1+t_1)}N_R(m-t-m_1+t_1,n-m_1+t_1)N_R(t-t_1,k-t_1).$$
\end{thm}
\proof
By counting the number of couples $(P,Q)$, where $P\in{\cal M}_R(m_1,t_1;n+k,n)$ and $Q\in{\cal M}_R(m,t;n+k,n)$
  with $P\subseteq Q$, by Theorems~\ref{lem5.1} and~\ref{lem5.2},
  $$N'_R(m_1,t_1;m,t;n+k,n)=\frac{N_R(m,t;n+k,n)N_R(m_1,t_1;m,t;n+k,n)}{N_R(m_1,t_1;n+k,n)}.$$
Therefore, the desired result follows. $\qed$

 \section{Arcs and caps}
 In this section, we always assume that $R=R_1\times R_2\times\cdots\times R_\ell$ is a commutative ring,
where $R_1,R_2,\ldots,R_\ell$ are local rings  with maximal ideals $M_1, M_2,\ldots, M_\ell$, respectively.
Recall that $\rho_i$ is the projection map from $R$ to $R_i$, and
let  $\pi_i$ be the natural surjective homomorphism from
$R_i$  to $R_i/M_i$ for all $i\in[\ell]$.
Recall that $R^{n}$ is the $n$-dimensional row vector space over $R$.
The $1$-subspaces, $2$-subspaces, and $(n-1)$-subspaces of $R^{n}$ are called {\it points, lines},
and {\it hyperplanes}  of $R^{n}$, respectively.

A $k$-{\it arc} in $R^n$ is a set ${\cal A}$ of $k$ points with the property that any $n$ of them
span the whole space. An arc is {\it complete} if it cannot be extended to a larger arc.
An arc of size $\leq n$ in $R^n$ is always incomplete, and therefore we will
implicitly assume that an arc has size at least $n + 1$.
The maximum number of points that a $k$-arc can have is $m(n,R)$.
An arc in $R^3$ is called a {\it planar arc}.

A $k$-{\it cap} in $R^n$ is a set ${\cal C}$ of $k$ points  with the property that any $3$ of them span a $3$-subspace of $R^n$.
We always assume that $n\geq3$.
A cap is {\it complete} if it cannot be extended to a larger cap.
The maximum number of points that a $k$-cap can have is $m_2(n,R)$.
A cap in $R^3$ is called a {\it planar cap}. Planar caps and  planar arc are equivalent, and so $m(3,R)=m_2(3,R)$.

Arcs and caps in the vector
space $\mathbb{F}_q^n$ over the finite field $\mathbb{F}_q$ have a long history.
In 1947, Bose \cite{Bose} was the first to observe that the largest planar cap over $\mathbb{F}_q$ has size $q + 1$ if
$q$ is odd and $q + 2$ if $q$ is even. In 1955-1957,  Segre \cite{Segre1}-\cite{Segre3}
published his fundamental work on arcs. For more results about arcs and caps, see \cite{Hirschfeld,Hirschfeld2}.
The most fundamental problem  is that of determining the size and
the structure of the largest arcs (resp. caps) in $\mathbb{F}_q^n$.

Let ${\cal S}\subseteq{\cal M}_R(1,n)$, where ${\cal M}_R(1,n)$ is the set of all the points of $R^n$. Define
$\rho_i({\cal S})=\{\langle\rho_i(\alpha)\rangle : \langle\alpha\rangle\in{\cal S}\}$ for all $i\in[\ell]$,
and $\pi_i(\rho_i({\cal S}))=\{\langle\pi_i(\rho_i(\alpha))\rangle : \langle\alpha\rangle\in{\cal S}\}$ for all $i\in[\ell]$.

\begin{thm}\label{lemma6.1}
Let $n\geq 2$ and ${\cal A}$ be an arc in $R^n$. Then
\begin{itemize}
\item[\rm(i)] $\rho_i({\cal A})$ is an $|{\cal A}|$-arc in $R_i^n$,
and $\pi_i(\rho_i({\cal A}))$ is an $|{\cal A}|$-arc in $(R_i/M_i)^n$  for all $i\in[\ell]$.

\item[\rm(ii)]
${\cal A}$ is complete if and only if there exists some $i_0\in[\ell]$ such that $\pi_{i_0}(\rho_{i_0}({\cal A}))$ is complete.
\end{itemize}
\end{thm}
\proof (i) Since ${\cal A}$ is an arc in $R^n$, by Lemma~\ref{lemma3.1}, $\rho_i({\cal A})$ is an arc  in $R_i^n$ and $|{\cal A}|=|\rho_i({\cal A})|$  for all $i\in[\ell]$. By Corollary~\ref{lemma2.1-2},
$\pi_i(\rho_i({\cal A}))$ is an arc  in $(R_i/M_i)^n$ and $|{\cal A}|=|\pi_i(\rho_i({\cal A}))|$ for all $i\in[\ell]$.

(ii) Suppose that there exists some $i_0\in[\ell]$ such that $\rho_{i_0}({\cal A})$ is complete. If ${\cal A}$ is not complete, then there exists
some $\langle\alpha\rangle\in{\cal M}_R(1,n)$ such that ${\cal A}\cup\{\langle\alpha\rangle\}$ is an arc in $R^n$. By (i),
$\rho_{i_0}({\cal A})\cup\{\langle\rho_{i_0}(\alpha)\rangle\}$ is an arc in $R_{i_0}^n$, a contradiction since $\rho_{i_0}({\cal A})$ is complete.
Conversely, suppose that ${\cal A}$ is complete. If $\rho_i({\cal A})$ is not complete  for all $i\in[\ell]$, then there is an
$\langle \alpha_i\rangle\in {\cal M}_{R_i}(1,n)$ such that $\rho_i({\cal A})\cup\{\langle\alpha_i\rangle\}$ is an arc in $R_i^n$.
From (\ref{shi3}), we deduce that there exists the unique $\langle\alpha\rangle\in{\cal M}_R(1,n)$ such that $\langle\rho_i(\alpha)\rangle=\langle\alpha_i\rangle$ for all $i\in[\ell]$. By Lemma~\ref{lemma3.1}, ${\cal A}\cup\{\langle\alpha\rangle\}$
is an arc in $R^n$, a contradiction since ${\cal A}$ is complete.

By Corollary~\ref{lemma2.1-2}, $\rho_i({\cal A})$ is complete if and only if $\pi_i(\rho_i({\cal A}))$ is complete  for all $i\in[\ell]$.
Therefore, the desired result follows. $\qed$

\begin{thm}\label{lemma6.2}
 $m(n,R)=\min\{m(n,R_i/M_i): i\in[\ell]\}$.
 \end{thm}
\proof Suppose that  ${\cal A}$ is an $m(n,R)$-arc in $R^n$. By Theorem~\ref{lemma6.1}, $\rho_i({\cal A})$
is an $m(n,R)$-arc in $R_i^n$ for all $i\in[\ell]$, which implies that $m(n,R)\leq m(n,R_i)$ for all $i\in[\ell]$.
Conversely, if $m=\min\{m(n,R_i): i\in[\ell]\}$,
then there exists an   $m$-arc ${\cal A}_i$ in $R_i^n$ such that
$${\cal A}_i=\{\langle\alpha_{i1}\rangle,\alpha_{i2}\rangle,\ldots,\langle\alpha_{im}\rangle\}\quad\hbox{for all}\;i\in[\ell].$$
From (\ref{shi3}), we deduce that there exists the unique $\langle\alpha_j\rangle\in{\cal M}_R(1,n)$ such that $\langle\rho_i(\alpha_j)\rangle=\langle\alpha_{ij}\rangle$ for all $i\in[\ell]$ and $j\in[m]$. By Lemma~\ref{lemma3.1},
$\{\langle\alpha_1\rangle,\langle\alpha_2\rangle,\ldots,\langle\alpha_m\rangle\}$ is an $m$-arc in $R^n$, which implies that
$m\leq m(n,R)$. So, $m(n,R)=\min\{m(n,R_i): i\in[\ell]\}$.

By Corollary~\ref{lemma2.1-2},  $m(n,R_i)=m(n,R_i/M_i)$ for all $i\in[\ell]$.  Therefore, the desired result follows.
  $\qed$

From \cite[Section~6]{Hirschfeld2},   some values $m(n,\mathbb{F}_q)$\footnote{Note that $m(n,\mathbb{F}_q)$ is denoted by $m(n-1,q)$ in \cite{Hirschfeld2}.} in $\mathbb{F}_q^n$ obtained are as follows:
\begin{eqnarray*}
m(2,\mathbb{F}_q)&=&q+1.\\
m(3,\mathbb{F}_q)&=&q+\frac{3+(-1)^q}{2}.\\
m(4,\mathbb{F}_q)&=&q+1\;\hbox{for}\;q>3.\\
m(5,\mathbb{F}_q)&=&q+1\;\hbox{for $q\geq5$}.\\
m(n+1,\mathbb{F}_q)&=&n+2\;\hbox{for}\;q\leq n+1.\\
\end{eqnarray*}

By Theorem~\ref{lemma6.2} and above the results, we obtain the following corollary.

\begin{cor}\label{lemma6.3}
Let $n\geq2$ and $q_i=|R_i|/|M_i|$ for all $i\in[\ell]$. Then
\begin{eqnarray*}
m(2,R)&=&\min\{q_i+1:i\in[\ell]\}.\\
m(3,R)&=&\min\left\{q_i+\frac{3+(-1)^{q_i}}{2} : i\in[\ell]\right\}.\\
m(4,R)&=&\min\{q_i+1:i\in[\ell]\}\quad\hbox{if}\;q_i> 3\;\hbox{for all}\;i\in[\ell].\\
m(5,R)&=&\min\{q_i+1:i\in[\ell]\}\quad\hbox{if}\;q_i\geq 5\;\hbox{for all}\;i\in[\ell].\\
m(n+1,R)&=&n+2\quad\hbox{if}\;q_i\leq n+1\;\hbox{for all}\;i\in[\ell].
\end{eqnarray*}
\end{cor}

\begin{thm}\label{lemma6.4}
Let $n\geq 3$ and ${\cal C}$ be a cap in $R^n$. Then
\begin{itemize}
\item[\rm(i)] $\rho_i({\cal C})$ is a $|{\cal C}|$-cap in $R_i^n$,
and $\pi_i(\rho_i({\cal C}))$ is a $|{\cal C}|$-cap in $(R_i/M_i)^n$  for all $i\in[\ell]$.

\item[\rm(ii)]
${\cal C}$ is complete if and only if there exists some $i_0\in[\ell]$ such that $\pi_{i_0}(\rho_{i_0}({\cal C}))$ is complete.
\end{itemize}
\end{thm}
\proof The proof is similar to that of Theorem~\ref{lemma6.1}, and will be omitted. $\qed$

\begin{thm}\label{lemma6.5}
 $m_2(n,R)=\min\{m_2(n,R_i/M_i): i\in[\ell]\}$.
 \end{thm}
 \proof The proof is similar to that of Theorem~\ref{lemma6.2}, and will be omitted. $\qed$

From \cite[Section~6]{Hirschfeld2},   some values $m_2(n,\mathbb{F}_q)$\footnote{Note that $m_2(n,\mathbb{F}_q)$ is denoted by $m_2(n-1,q)$ in \cite{Hirschfeld2}.} in $\mathbb{F}_q^n$ obtained are as follows:
\begin{eqnarray*}
m_2(n+1,\mathbb{F}_2)&=&2^n.\\
m_2(3,\mathbb{F}_q)&=&q+\frac{3+(-1)^q}{2}.\\
m_2(4,\mathbb{F}_q)&=&q^2+1\;\hbox{for}\;q>2.\\
m_2(5,\mathbb{F}_3)&=&20.\\
m_2(6,\mathbb{F}_3)&=&56.\\
m_2(5,\mathbb{F}_4)&=&41.
\end{eqnarray*}

By Theorem~\ref{lemma6.5} and above the results, we obtain the following corollary.

\begin{cor}\label{lemma6.6}
Let $n\geq3$ and $q_i=|R_i|/|M_i|$ for all $i\in[\ell]$. Then
\begin{eqnarray*}
m_2(n+1,R)&=&2^{n}\quad\hbox{if}\;q_i= 2\;\hbox{for all}\;i\in[\ell].\\
m_2(3,R)&=&\min\left\{q_i+\frac{3+(-1)^{q_i}}{2} : i\in[\ell]\right\}.\\
m_2(4,R)&=&\min\{q_i^2+1; i\in[\ell]\}\quad\hbox{if}\;q_i> 2\;\hbox{for all}\;i\in[\ell].\\
m_2(5,R)&=&20\quad\hbox{if}\;q_i=3\;\hbox{for all}\;i\in[\ell].\\
m_2(6,R)&=&56\quad\hbox{if}\;q_i=3\;\hbox{for all}\;i\in[\ell].\\
m_2(5,R)&=&41\quad\hbox{if}\;q_i=4\;\hbox{for all}\;i\in[\ell].
\end{eqnarray*}
\end{cor}

\section{Concluding remarks}
Recall that $R=R_1\times R_2\times\cdots\times R_\ell$ is a commutative ring and $\rho_i$ is the projection map from $R$ to $R_i$ for all $i\in[\ell]$, where $R_i$ is a local ring  with the unique maximal ideal $M_i$.

For symplectic spaces and orthogonal spaces
of odd characteristic, Wan \cite{Wan2} obtained Anzahl formulas of  all subspaces   over the finite field $\mathbb{F}_q$;
Sirisuk and Meemark \cite{Sirisuk} obtained Anzahl formulas of totally isotropic subspaces over $R$.  In \cite{Guo1,Guo2}, we
generalized Sirisuk-Meemark's results  over $R$.

For $0_{m,n} \not= A \in R^{m\times n}$, by Cohn's definition \cite{Cohn},
the {\it inner rank} of $A$, denoted by ${\rm ir}(A)$,
is the least positive integer $r$ such that $A = BC$ where $B \in R^{m\times r}$ and $C\in R^{r\times n}$.
Define ${\rm ir}(0_{m,n})= 0$. By \cite[Proposition~2.1]{Sirisuk2}, ${\rm ir}(A)=\max\{{\rm ir}(\rho_i(A)):1\leq i\leq \ell\}$.
For the $n$-dimensional row vector space $\mathbb{Z}_h^n$ over the residue class ring $\mathbb{Z}_h$,
by \cite[Theorem~3.3]{Guo3}, we obtain the following dimensional formula
$$\dim(P) + \dim(Q) - \dim(P\cap Q) = {\rm ir}\left(\begin{array}{c}P\\ Q\end{array}\right),$$
where $P$ and $Q$ are two subspaces of $\mathbb{Z}_h^n$. It seems interesting to discuss the dimensional formula
  in the $n$-dimensional row vector space $R^n$.

Let $0\leq r \leq m \leq n$.
A family ${\cal F}_R\subseteq{\cal M}_R(m,n)$ is called $r$-{\it intersecting} if $\dim(A\cap B)\geq r$ for all $A, B\in{\cal F}_R$.
Frankl and Wilson \cite{Frankl}  obtained an upper bound
on the size of an $r$-intersecting family
and described exactly which families meet this bound over the finite field $\mathbb{F}_q$. Huang, Lv and Wang  \cite{Huang3} and the first present author \cite{Guo3} discussed similar problems over the residue class ring $\mathbb{Z}_{h}$. It seems interesting to discuss the $r$-intersecting families  over the commutative ring $R$.

\section*{Acknowledgments}
This research is supported by
 National Natural Science Foundation of China (11971146)
 and the Fundamental Research Funds for the Universities in Hebei Province (JYT202102).

\end{document}